\documentclass{amsart}
\usepackage{amssymb}

\newcommand{\XREF}{\label}

\newcommand{\R}{{\mathbb R}}
\newcommand{\C}{{\mathbb C}}

\newcommand{\CC}{{\overline \C}}

\newcommand{\A}{{\mathcal A}}
\newcommand{\downto}{{\searrow}}

\newcommand{\nbhd}{{neighborhood }}
\newcommand{\tri}{{\bigtriangleup}}

\newtheorem{corollary}{Corollary}
\newtheorem{lemma}{Lemma}

\newtheorem{theorem}{Theorem}
\newtheorem{proposition}{Proposition}

\newtheorem{othertheorem}{Theorem}

\theoremstyle{definition}
\newtheorem{definition}{Definition}
\theoremstyle{remark}
\newtheorem{remark}{Remark}

\newtheorem{example}{Example}

\numberwithin{equation}{section}
\numberwithin{theorem}{section}
\numberwithin{definition}{section}
\numberwithin{remark}{section}
\numberwithin{example}{section}
\numberwithin{lemma}{section}
\numberwithin{property}{section}
\numberwithin{proposition}{section}
\numberwithin{claim}{section}
\numberwithin{othertheorem}{section}
\numberwithin{conj}{section}
\numberwithin{corollary}{section}

\begin{document}

\title{Uniformly perfect sets, Rational
semigroups, Kleinian groups and IFS's}
\author{Rich Stankewitz}
\address{Department of Mathematics, 
         Texas A\&M University, 
         College Station, TX 77843}
\email{richs@math.tamu.edu}

\begin{abstract}
We show that the Julia set of a non-elementary rational semigroup $G$
is uniformly perfect when there is a uniform bound on the Lipschitz
constants of the generators of $G$.  This also proves that the limit
set of a non-elementary M\"obius group is uniformly perfect when there
is a uniform bound on the Lipschitz constants of the generators of the
group and this implies that the limit set of a finitely generated
non-elementary Kleinian group is uniformly perfect.
\end{abstract}

\maketitle

\section{Introduction}
A rational semigroup $G$ is a semigroup 
of nonconstant rational functions 
defined on the Riemann sphere $\CC$ with the semigroup operation being
functional composition.  When a semigroup $G$ is generated by the functions
$\{f_i:i \in I\}$, for some index set $I$, 
we write this as
\begin{equation}
G=\langle f_i:i \in I \rangle.\notag
\end{equation}
Note that in~\cite{HM1} and~\cite{HMJ} rational semigroups are always
taken to have at least one element of degree at least two.  We do not
make such a restriction here.

The study of rational semigroups may be viewed as a generalization of
the study of Kleinian groups, iteration of a rational function, and
iterated function systems (IFS's). 
For example limit sets of Kleinian groups,
Julia sets of a rational function, and self-similar sets generated by
iterated function systems are Julia sets of rational semigroups
(see Remarks~\ref{rmkKG} and~\ref{rmkIFS}) and in
these sets fixed points are always dense (see~\cite{Be3},
\cite{Be}, and~\cite{Hu}). 

It is known that the Julia set of a rational function is uniformly 
perfect. Several proofs of this fact have been given, namely by
Eremenko~\cite{Er}, Hinkkanen~\cite{Hi}, and Ma\~n\'e and da Rocha~\cite{MdR}.
Also, Hinkkanen and Martin have shown in~\cite{HMJ} that Julia sets of
finitely generated rational semigroups are uniformly perfect when the
semigroup is generated by maps all which have degree greater than or
equal to two.  (Although in~\cite{HMJ} the theorem is stated
under the weaker hypothesis that only one of the maps of the semigroup
must have degree at least two, a closer inspection shows that the
proof is only valid when all the maps are of degree at least two.)  
The proof uses the fact that Julia sets of single rational
functions of degree at least two 
are uniformly perfect and also relies substantially on the
fact that the semigroup is finitely generated.
In this paper we show in Theorem~\ref{main}
that when there exists a uniform bound on the Lipschitz
constants with respect to the spherical metric associated with each generator
of the rational semigroup the corresponding Julia set is
necessarily uniformly perfect.  This Lipschitz condition is trivially
satisfied when the semigroup is finitely generated and so we see that
this is an improvement on the previous result.  Furthermore the 
proof does not rely
on the fact that Julia sets of single rational functions are uniformly
perfect and we also relax the restriction on the degree of the maps in the
semigroup, thus we are able to apply this
result to M\"obius groups as well.

The author would like to thank the referee for the helpful comments.

\section{Definitions and basic facts}
In~\cite{HM1}, p.~360 the definitions of the set of
normality, often called the Fatou set, 
and the Julia set of a rational semigroup are as follows:

\begin{definition} \label{N(G),J(G)}
For a rational semigroup $G$ we
define the set of normality of $G$, $N(G)$, by 
$$N(G)=\{z \in \CC:\text{there is a \nbhd of } z \text{ on
which } G \text{ is a normal family} \}$$
and define the Julia set of $G$, $J(G)$, by
$$J(G)=\CC \setminus N(G).$$
\end{definition}

When $G=\langle f \rangle$, we abuse the notation and write $N(f)$ and
$J(f)$ for $N(\langle f \rangle)$ and $J(\langle f \rangle)$
respectively.

\begin{proposition}[\cite{HM1}, p.~360]\label{invprop}
The set $N(G)$ is forward invariant under
each element of $G$ and $J(G)$ is backward invariant under each
element of $G$.
\end{proposition}

The following proposition about the backward self-similarity of the
Julia set has been noted by many people and is easy to prove.

\begin{proposition}\label{backselfsim}
If $G=\langle g_1, \dots, g_N \rangle$, then $J(G)=\cup_{i=1}^N
g_i^{-1}(J(G))$ and $N(G)=\cap_{i=1}^N g_i^{-1}(N(G))$.
\end{proposition}

The sets $N(G)$ and $J(G)$ are, however, not necessarily 
completely invariant under the elements of $G$.  
This is in contrast to the case of single function dynamics.  For a
discussion on completely invariant Julia sets the reader is referred
to~\cite{RS},~\cite{RS1} and~\cite{RS2}.   

\begin{definition}\label{non-el}
We say that a rational semigroup $G$ is non-elementary if $J(G)$ has
three or more points.
\end{definition}

\begin{lemma}[\cite{HM1}, Lemma 3.1]\label{J(G)perfect}
The set $J(G)$ is perfect when $G$ is non-elementary.
\end{lemma}

Using an idea of Baker in~\cite{B}, Hinkkanen and Martin have shown
the following theorem.  (The author has shown that the idea of Schwick
in~\cite{Sch} can also be adapted to prove the following result.)

\begin{othertheorem}[\cite{HM1}, Theorem~3.1 and Corollary~3.1] \label{T:repdense}
If $G$ is a non-elementary 
rational semigroup, then the repelling fixed points of the
elements of $G$ are dense in $J(G)$.  Hence also
$$J(G)= \overline{\bigcup_{g \in G} J(g)}.$$
\end{othertheorem}

\begin{remark}
In~\cite{HM1} 
all the semigroups $G$ are assumed to have at least one
element of degree two or more, but Lemma~\ref{J(G)perfect} and
Theorem~\ref{T:repdense} still hold
(with the proofs given in the references) if we remove this restriction.
\end{remark}

\begin{remark}\label{rmkKG}
A M\"obius group $G$, i.e., a group of maps of the form
$(az+b)/(cz+d)$ where $ad-bc \neq 0$, is a rational semigroup.  If
there does not exist a finite $G$-orbit in $\R^3$ (this is the
definition of non-elementary in~\cite{Be3}) then one can see
by~\cite{Be3}, p.~90 that $G$ is non-elementary in the sense of
Definition~\ref{non-el}.  Whereas, if a M\"obius group $G$ is
non-elementary in the sense of Definition~\ref{non-el}, then $G$ is
non-elementary in the sense of~\cite{Be3} or $G$ is conjugate to a
group of linear transformations (see~\cite{Be3}, p. 84).  

Since the limit set of a non-elementary M\"obius group (which need not
be discrete) is the closure of the set of repelling fixed points
(see~\cite{Be3}, p.~97) we see by Theorem~\ref{T:repdense} that the
notion of limit set and Julia set coincide in the case of
non-elementary M\"obius groups and, in particular, for non-elementary
Kleinian groups, i.e., M\"obius groups whose action is properly discontinuous
at some point $z \in \CC$.  .
\end{remark}

\begin{remark}\label{rmkIFS}
An iterated function system (IFS) generated by the (possibly infinite)
set of contracting linear maps $\{g_i: i \in I
\}$ on $\C$ (see~\cite{Hu}) corresponds to a rational semigroup
generated by the inverses of the generating maps $\{g_i^{-1}:i \in I
\}$.  Note that the attracting fixed point of an element $g$ in the
IFS is a repelling fixed point for $g^{-1}$ in the corresponding
semigroup.  Defining the attractor set $A$ for an IFS to be the
closure of the attracting fixed points (see~\cite{Hu}), we see that
$A$ is the Julia set for the corresponding rational semigroup.  One
may also note that when the IFS is generated by $\{g_1, \dots, g_n\}$
the attractor set $A$ may be defined as the unique compact set $A
\subset \C$ which satisfies $A=\cup_{i=1}^N g_i(A)$ and so it also
follows from Proposition~\ref{backselfsim} that $A$ is the Julia set
for the corresponding rational semigroup.
\end{remark}

\begin{example}
Define the maps $g_1(z)=z/3$, $g_2(z)=(e^{\pi i/3}z+1)/3$, $g_3(z)=(e^{2 \pi
i/3}z+2)/3$, 
and $g_4(z)=(z+2)/3.$  Then $J(\langle g_1^{-1}, g_4^{-1} \rangle)$
is the middle third Cantor set and 
$J(\langle g_1^{-1}, g_2^{-1}, g_3^{-1}, g_4^{-1}\rangle)$ 
is the von Koch curve (see~\cite{F}, p.~xv).
\end{example}

\begin{definition}
For a rational function $g$ of degree two or more and a point $z \in
\CC$ we define the backward orbit with respect to $g$ by 
$O_g^-(z)=\{w: \text{there exists }  n  \text{ such that }
g^n(w) = z\}$.
\end{definition}

\begin{definition}
If $g$ is a rational function of degree two or more we
define the exceptional set of $g$ to be $E(g)=\{z\in \CC:O_g^-(z)
\text{ contains at most two points} \}.$
If $g$ is a M\"obius map with an attracting fixed point 
and a repelling fixed point, i.e., a
loxodromic M\"obius map, then we
define the exceptional set $E(g)$ of $g$ to be the set consisting of
the single attracting fixed point.
\end{definition}

\begin{proposition}[\cite{Be}, Theorem 6.9.4] \label{P:singexp}
        Let $f$ be a rational function with $\deg f \ge 2$ or a
        loxodromic M\"obius map.  Let $W$ be a 
        non-empty open set intersecting $J(f)$, and let $K$ be a compact 
        subset of $\CC \setminus E(f)$.  Then there exists an integer $N$ 
        such that $K \subset f^n(W)$ for all $n \ge N$.
\end{proposition}

\begin{definition}
For a rational function $g$ we define 
Lip $g = \inf \{M:\sigma(g(z),g(w))\leq M
\sigma(z,w)$ for all $z,w \in \CC\}$ where $\sigma$ denotes the
spherical metric.
\end{definition}

Note that for a rational map $g$, we have Lip $g < +\infty$ as can be
seen by proving the equivalent statement that $\frac
{|g'(z)|(1+|z|^2)}{1+|g(z)|^2}$ is bounded above on $\CC$
(see~\cite{Be}, p.~32).

\begin{definition} \XREF{confann}
A conformal annulus is an open subset $\A$ of $\CC$ that can be conformally
mapped onto the genuine annulus $Ann(0;r_1,r_2)=\{z: 0 \leq r_1 < |z| < 
r_2 \leq \infty\}$
and the modulus of such a conformal annulus is given by
$$\text{mod}(\A)=\frac{1}{2\pi} \log\frac{r_2 }{ r_1}.$$
\end{definition}

We note that $\text{mod}(\A)$ is a conformal invariant. 

\begin{definition} \XREF{sepann}
A conformal annulus $\A$ is said to separate a set $F$ if $F$ intersects
both components of $\CC \setminus \A$ and $F \cap \A = \emptyset$.
\end{definition}

\begin{definition} \XREF{updef}(\cite{P}, p.~192)
We say that a compact subset $F \subset \CC$ is uniformly perfect if
$F$ has at least two points and if there is a uniform upper bound on the 
moduli of all conformal annuli in
$\CC \setminus F$ which separate $F$.  
\end{definition}

Uniformly perfect sets
were introduced by A. F. Beardon and Ch.~Pommerenke in 1978 in \cite{BP}.

\section{The main result}

\begin{theorem}\label{main}
Let $G=\langle g_i:i \in I \rangle$ be a non-elementary
rational semigroup generated by the maps $\{g_i:i \in I\}$ 
such that $\sup_{i \in I}$ Lip $g_i \leq C<+\infty$. 
Then the Julia set $J(G)$ is uniformly perfect.
\end{theorem}

\begin{proof}
The proof given below is an adaptation of the proof that the Julia set
of a single rational function of degree at least two is uniformly
perfect which was presented in L.~Carleson and T.~Gamelin's
book~\cite{CG}.  In this proof all notions of distance and convergence
will be with respect to the spherical metric.  Fix four points $z_1,
z_2, z_3$ and $z_4$ in $J(G)$. Let $\delta >0$ be chosen small enough
so that any two of the four selected points is at a (spherical)
distance strictly greater than $C \delta$ from each other and such
that $\delta < \frac{4d}{5}$ where $d=diam(J(G))$ denotes the
spherical diameter of $J(G)$.

Suppose there is a sequence $\{A_n\}$ of conformal annuli in $N(G)$
with moduli tending to $\infty$ such that both components of $\CC
\setminus A_n$ meet $J(G)$.  Let $E_n$ be the component of $\CC
\setminus A_n$ with the smaller spherical diameter.  The diameter
$diam(E_n)$ tends to zero as can be seen by noting that that there are
simple closed curves in $A_n$ whose hyperbolic lengths tend to 0 and
which separate $E_n$ from $\CC \setminus (A_n \cup E_n)$.  We assume
that all $diam(E_n) < \delta$.

Since $A_n \cup E_n$ is open and meets $J(G)$, we know by
Theorem~\ref{T:repdense} that there exists a point $z \in J(G)
\cap (A_n \cup E_n)$ which is a repelling fixed point for the map $h
\in G$, say.  If $E(h) \cap J(G) = \emptyset$, where $E(h)$ denotes the set of 
(at most two) exceptional points of $h$, let $U = \emptyset$.
Otherwise, let $U$ be a union of at most two open spherical discs
centered at the points of $E(h) \cap J(G)$ each of diameter at most
$\frac{d}{10}$ such that the boundary of each disk in $U$ contains a
point of $J(G)$.  Such disks can be found since $J(G)$ is perfect by
Lemma~\ref{J(G)perfect}.  Let $K=J(G) \setminus U$ and note that
$diam(K) \geq \frac{4d}{5}>\delta$.  Then by the expanding property of
Julia sets (Proposition~\ref{P:singexp}) there exists a positive
integer $k$ such that $h^k(A_n \cup E_n) \supset K$.  Since $h^k(A_n)
\subset N(G)$, which follows directly from Proposition~\ref{invprop},
we conclude $h^k(E_n) \supset K$ and so $diam(h^k(E_n)) > \delta$.  We
now let $f_n \in G$ be a function of minimal word length such that
$diam(f_n(E_n)) > \delta$.  Each $f_n$ can be written in the form
$$f_n = g_{i_1} \circ g_{i_2} \circ \dots \circ g_{i_m}$$ where each
$i_j \in I$ depends on $n$.  Letting $F_n = g_{i_2} \circ \dots \circ
g_{i_m}$ (we let $F_n = Id$ if $f_n$ is one of the generating
functions in the set $\{g_i:i \in I\}$) we see by the minimality of
the choice of $f_n$ that $diam(F_n(E_n)) \leq
\delta$.  Hence $diam(f_n(E_n)) =diam(g_{i_1}(F_n(E_n))) \leq 
C diam(F_n(E_n)) \leq C \delta.$  

Let $\phi_n$ be a conformal map from the unit disk $\tri$ onto $A_n
\cup E_n$ such that $\phi_n(0) \in E_n$.  Note that here we use the fact that
$J(G)$ is perfect to assert that $\CC \setminus (A_n \cup E_n)$
contains more than two points and therefore $A_n \cup E_n$ is
conformally equivalent to the unit disk.  Let $K_n=\phi_n^{-1} (E_n)$
and note that as the modulus of $\tri \setminus K_n$ equals the
modulus of $A_n$, we conclude as above that $diam(K_n) \to 0$.  Let
$h_n=f_n \circ \phi_n$.  Note that $h_n(\tri \setminus K_n) \subset
N(G)$ and $diam(h_n(K_n))=diam(f_n(E_n)) \leq C \delta$.  Hence each
$h_n(\tri )$ omits at least three of the points $z_1, z_2, z_3$ and
$z_4$ in $J(G)$.  Hence $\{h_n\}$ is normal in $\tri$ by Montel's
theorem.  By the equicontinuity of the family $\{h_n\}$ we have
$diam(h_n(K_n)) \to 0$ since $diam(K_n) \to 0$ and $0 \in K_n$ for all
$n$.  This is a contradiction since $diam(h_n(K_n))=diam(f_n(E_n)) >
\delta$.
\end{proof}


Let $G'$ be a M\"obius group generated (as a group) by $\{g_i:i \in
I\}$.  Since the semigroup generated by $\{g_i: i \in I\} \cup
\{g_i^{-1}:i \in I\}$ is $G'$ and Lip $g_i=$ Lip $g_i^{-1}$
when $g_i$ is a M\"obius transformation (see~\cite{Be}, p.~33) we see that
Theorem~\ref{main} implies the following corollary.

\begin{corollary}
The limit set (Julia set) 
$J(G')$ of a non-elementary M\'obius group is uniformly
perfect when $\{g_i:i \in
I\}$, the generators of $G'$ (as a group), satisfy 
$\sup_{i \in I}$ Lip $g_i<+\infty$.  
\end{corollary}

We note that Kleinian groups which possess a uniform bound on the
Lipschitz constants of its generators must necessarily be finitely
generated since otherwise one can show that the group is not discrete
(see~\cite{Be}, p. 33 or~\cite{Be3}, p. 42).  Hence Theorem~\ref{main}
can only be used to duplicate the following result due to Pommerenke
(see~\cite{P2}).  Pommerenke's proof is based on more analytic
methods as opposed to the more geometric view used in this paper.

\begin{corollary}
The limit set (Julia set) 
$J(G')$ of a non-elementary Kleinian group is uniformly
perfect when $G'$ is finitely generated.
\end{corollary}

\begin{remark}
Not all limit sets of Kleinian groups are uniformly perfect.
For each positive integer $n$ 
let $C_n =\{z:|z-a_n|=r_n\}$ and $C_n'=\{ z:|z-a_n -2i|=r_n \}$ where the real 
numbers $a_n \downto 0$, $r_n \downto 0$ such that $r_n <a_n \leq 1$ and 
$\frac{a_n-r_n}{a_{n+1}+r_{n+1}} \to +\infty$.  Then one can show
that the Schottky group 
generated by the M\"obius maps $g_n$ which map $C_n$ onto $C_n'$ taking 
the interior of $C_n$ onto the exterior of $C_n'$ has a limit set which 
is not uniformly perfect.  One can see this by noting that the 
annuli $A_n=\{z:a_{n+1}+r_{n+1} < |z| < a_n -r_n\}$ separate the limit set.
\end{remark}

\begin{corollary}
The attractor set $A$ of an IFS generated by the contracting linear
maps $\{g_i:i \in I\}$ defined on $\C$ is uniformly perfect when
$\sup_{i \in I}$ Lip $g_i<+\infty$.
\end{corollary}

\section{Applications}

\begin{othertheorem}[\cite{HMJ}, Theorem~4.1]\label{sainJ(G)}
Let $G$ be a rational semigroup such that $J(G)$ is uniformly perfect.
Suppose that $z_0$ is a superattracting fixed point of an element $h
\in G$.  Let $A$ be the union of all the components of $N(h)$ in which
the iterates of $h$ tend to $z_0$.  Then either $z_0
\in N(G)$ or $A \subset J(G).$  In particular, either $z_0
\in N(G)$ or $z_0$ lies in the interior of $J(G)$.
\end{othertheorem}

\begin{corollary}
Let $G=\langle g_i:i \in I \rangle$ be a non-elementary
rational semigroup generated by the maps $\{g_i:i \in I\}$ 
such that $\sup_{i \in I}$ Lip $g_i <+\infty$ and let $z_0$ be a
superattracting fixed point of some element of $G$.  Then either $z_0$
lies in (the interior of) the Fatou set of $G$ or in the interior of
the Julia set of $G$.
\end{corollary}

Julia sets of rational semigroups which contain elements of degree two or 
more are, however, not always uniformly perfect.  In fact
the following theorem is true.

\begin{othertheorem}[\cite{HMJ}, Theorem 5.1]
There exists an infinitely generated rational semigroup $G$ (all of
whose elements have degree at least two) with the property that for
any positive integer $N$, the semigroup $G$ contains only finitely
many elements of degree at most $N$, such that $J(G)$ is not uniformly
perfect, and such that $G$ contains an element $g$ with a
superattracting fixed point $\alpha$ with $\alpha \in \partial J(G)
\subset J(G).$
\end{othertheorem}

\begin{example}\label{ex}
Let $f_0$ be a non-constant rational function which has an attracting
or superattracting fixed point at $\infty$.  Pick a point $a \in
J(f_0)$ and let $b_n$ be a sequence of points in $\C$ tending to $a$.
Letting $f_n(z)=f_0(z+b_n)-b_n$, we can show that for the semigroup
$G=\langle z^2, f_n:n \geq 0 \rangle$ we have $\sup_{n \geq 0}$ Lip
$f_n <+\infty$ and consequently $\{z:|z| \leq 1\} \subset J(G)
\subsetneq \CC$.
\end{example}

\begin{proof}
Since the $b_n$'s are bounded one can easily show that there is a
small \nbhd of $\infty$ which maps into itself by every element of $G$ and
thus such a \nbhd must lie in $N(G)$.  Also since the $b_n$'s are
bounded one can show that $\sup_{n \geq 0}$ Lip $f_n <+\infty$ and so by
Theorem~\ref{main} $J(G)$ is uniformly perfect.  Since $a-b_n \in
J(f_n)$ for each $n$ we conclude that $0 \in
\overline{\cup_{n=0}^\infty J(f_n)} \subset J(G)$ and so by
Theorem~\ref{sainJ(G)} we have $\{z:|z| \leq
1\} \subset J(G)$.  
\end{proof}


One may ask whether there exists a finitely generated subsemigroup
$G'$ of $G$ in Example~\ref{ex} such that $0 \in J(G')$.  If so, one
could use the weaker result of Hinkkanen and Martin in~\cite{HMJ} (if
the degree of $f_0$ is greater than or equal to two) to conclude that
$J(G')$ is uniformly perfect and thus $\{z:|z| \leq 1\} \subset J(G')
\subset J(G)$.  It is often very difficult to accurately describe the
Julia set of a rational semigroup and so the above is a difficult
question.  The author would like to know if there exists an example of
a semigroup $G$ such that $J(G)$ contains a basin of attraction for
some superattracting fixed point of some element of $G$ such that
$J(G)
\subsetneq \CC$, but also such that no finitely generated subsemigroup
$G'$ is such that its Julia set contains this basin of attraction.

There are many open questions regarding the description of Julia sets of
rational semigroups that require further study in which theorems like
Theorem~\ref{main} and Theorem~\ref{sainJ(G)} may show to be useful.  
Some questions are the following:  What conditions imply that
$J(G)$ has
nonempty interior?  What conditions imply that 
$J(G)$ has empty interior?  When does
$J(G)$ have nonempty interior yet is a proper subset of $\CC$?
In~\cite{Su} and~\cite{Su2} Sumi uses an ``open set condition'' and
the backward self-similarity (Proposition~\ref{backselfsim}) to obtain 
some results in this direction and in particular shows that
certain conditions imply $J(G)$ may be a generalized Cantor set, 
have no interior, have zero
Lebesgue measure, or have Hausdorff dimension strictly less than 2.
See also~\cite{R}.

In order to get a feel for what certain Julia sets of finitely
generated rational semigroups may look like 
one may use computer algorithms to draw
these Julia sets by obtaining a measure whose support is exactly the Julia
set.  
Of course, the usual warnings must be heeded as certain situations exist
that seem to provide considerable barriers to the method giving an
accurate picture (see~\cite{Bo2}).

\bibliographystyle{amsplain}
\bibliography{upjset}

\end{document}